\pdfoutput=1

\documentclass{article}
\usepackage[T1]{fontenc}
\usepackage[UKenglish]{babel}

\usepackage{Factory, Math, Theorema, Classic, Styl}

\usepackage[all, 2cell]{xy} 
\UseAllTwocells 
\newcommand{\curvear}{\ar@/^1pc/} 		
\newcommand{\cruvear}{\ar@/^-1pc/}
\newcommand{\upar}{ \ar@<+.6ex> }		
\newcommand{\downar}{ \ar@<-.6ex> }
\newcommand{\dotar}{\ar@{.>}}
\newcommand{\updotar}{ \dotar@<+.6ex> }	
\newcommand{\downdotar}{ \dotar@<-.6ex> }

\newcommand{\con}[3]{\Big[\xymatrix@1{#1 & \ar[l]_-{#3} #2}\Big]}
\newcommand{\ext}[4]{\Big[\xymatrix@1{#1 \ar[r]^{#3}_-{#4} & #2}\Big]}
\newcommand{\lext}[4]{\Big[\xymatrix@1{#1 \ar[rr]^{#3}_-{#4} && #2}\Big]}
\newcommand{\maze}[6]{\Big[\xymatrix@1{#1 & \ar[l]_-{#4} #2 \ar[r]^-{#5}_-{#6} & #3}\Big]}
\newcommand{\lmaze}[6]{\Big[\xymatrix@1{#1 & \ar[l]_-{#4} #2 \ar[rr]^-{#5}_-{#6} && #3}\Big]}
\newcommand{\llmaze}[6]{\Big[\xymatrix@1{#1 & \ar[l]_-{#4} #2 \ar[rrr]^-{#5}_-{#6} &&& #3}\Big]}

\newcommand{\op}{\mathrm{op}}
\newcommand{\de}{\diamond} 			
\DeclareMathOperator*{\De}{\lozenge} 
\DeclareMathOperator{\ce}{\mathrm{cr}}

\newcommand{\Sets}{\mathfrak{Set}}
\newcommand{\Laby}{\mathfrak{Laby}}		

\newcommand{\FGrp}{\mathfrak{FGrp}}

\newcommand{\FCMon}{\mathfrak{FCMon}}
\newcommand{\FAb}{\mathfrak{FAb}}
\newcommand{\Ab}{\mathfrak{Ab}}
\newcommand{\FMod}{\mathfrak{FMod}}
\newcommand{\Mod}{\mathfrak{Mod}}

\DeclareMathOperator{\Add}{\mathrm{Add}}
\DeclareMathOperator{\Lin}{\mathrm{Lin}}
\DeclareMathOperator{\PMack}{\mathrm{PMack}}
\DeclareMathOperator{\Mack}{\mathrm{Mack}}
\DeclareMathOperator{\Pol}{\mathrm{Pol}}

\usepackage{mathtools}
\makeatletter
\newcommand\Label[1]{&\refstepcounter{equation}(\theequation)\ltx@label{#1}&}
\makeatother

\usepackage{fancyhdr}						
\xymatrixcolsep{1pc} 

\begin{document}


\titul{POLYNOMIAL FUNCTORS \\ ON POINTED CATEGORIES}
\auctor{Qimh Richey Xantcha\thanks{\textsc{Qimh Richey Xantcha}, Uppsala University: \texttt{qimh@math.uu.se}}}
\datum{\today}
\maketitle

\epigraph{
\begin{vverse}[Den, som vill skära glas, min herre, han måste hafva diamant!]
Den, som vill skära glas, min herre, han måste hafva diamant!

\vattr{Almqvist, \emph{Det går an}}
\end{vverse}
}

\bigskip 

\begin{abstract}
\noindent
A classification is provided of functors, in particular polynomial ones, 
from a category with a zero object in which every object is a finite sum of copies of a generating object, 
into an abelian category. 
This classification is extended to include functors 
from a category with sums and a zero object, carrying a small, 
regular projective generator. 

\end{abstract}

\bigskip

\noindent
Polynomial functors, first introduced by Eilenberg \& Mac Lane \cite{EM}, constitute today a rich theory 
with decisive applications to algebraic topology. Among the more significant developments of recent decades, 
one may mention the theory of unstable modules over the Steenrod algebra \cite{Steenrod} and the striking
applications to stable K-theory \cite{Scorichenko}.

Our purpose is, first, to provide a classification of functors $F\colon C\to D$, 
with $D$ being an abelian category and $C$ a 
\emph{pointed algebraic theory}: a 
category with a zero object in which every object is (isomorphic with) 
a finite sum of copies of a generating object $\Omega$. 

Results of this type abound in the scholarly literature. 
For example, Baues \& Pirashvili \cite{BP} found an equivalence between quadratic functors $\FGrp\to\Ab$ 
($\FGrp = $ finitely generated, free groups; $\Ab = $ abelian groups) and diagrams 
$$
\xymatrixcolsep{2pc}
\xymatrix{ M_e  \upar[r]^{H} & \upar[l]^{P} M_{ee}}
$$
of abelian groups $M_e$, $M_{ee}$ together with homomorphisms satisfying $PHP=2P$. Quadratic functors 
$\FAb\to\Ab$ ($\FAb = $ finitely generated, free abelian groups) 
are given by a diagrams of the same form, but where now $PHP=2P$ and $HPH=2H$ 
(Baues \cite{Baues}).

Baues, Dreckmann, Franjou \& Pirashvili \cite{BDFP} generalised this latter result to polynomial functors of 
any degree $n$, establishing an equivalence of categories 
$$
\Pol_n(\FAb,\Ab) \cong \Mack_n(\mathfrak{Sur},\Ab),
$$
the right-hand member denoting Mackey functors on
$\mathfrak{Sur}$, the category of finite sets and surjections, annihilated on sets of cardinality exceeding $n$.  
These Mackey functors essentially codify additive 
functors on spans of surjections; the exact mechanism is indicated hereunder. 

A peculiarity with the methods advocated in that paper, is that they do not bring about a tangible description 
of arbitrary functors $\FAb\to\Ab$, but instead proceed via functors $\FCMon\to\Ab$, with $\FCMon$ being 
the finitely generated, free, commutative \emph{monoids}, from which the polynomial functors $\FAb\to\Ab$ emerge as 
the ones annihilated on sets of cardinality exceeding $n$. 

The same phenomenon appears in the paper \cite{HPV}, by Hartl, Pirashvili \& Vespa, which derives an 
equivalence of categories 
$$
\Fun({}_P\mathfrak{FAlg}, \Ab) \cong \PMack(\Omega(P),\Ab),
$$
between functors on finitely generated, free algebras over a set operad $P$ and 
what they call \emph{pseudo-Mackey} functors on a certain subcategory $\Omega(P)$ of the 
May--Thomason category (\emph{ibid.},~Theorem 4.1). 
Again, polynomial functors correspond to those pseudo-Mackey functors annihilated on sets of cardinality 
exceeding $n$ (\emph{ibid.},~Corollary 4.9). 


An apposite classification of arbitrary functors $\FAb\to\Ab$ was not given until the work \cite{Laby}. 
In fact, letting $\Mod_A$ ($\FMod_A$) be the category of (finitely generated, free) right $A$-modules for some unital ring $A$, 
the theory expounded therein lays out an equivalence (Theorem 12)
$$
\Phi\colon \Fun(\FMod_A,\Mod_A) \to \Lin(\Laby(A),\Mod_A),
$$
encoding module functors as linear functors on the \emph{labyrinth category} $\Laby(A)$. 
The objects of $\Laby(A)$ are finite sets, and its arrows so-called \emph{mazes}.

Here is the principal idea. Letting $F\colon \FMod_A\to\Mod_A$ be a module functor,
one defines $\Phi(F)(X)$, for a set $X$, to be $\ce_X F(A,\dots,A)$, 
the cross-effect of $F$ of rank $X$ evaluated on $\abs{X}$ copies 
of $A$. The passages of a maze $P\colon Y\to X$ of $\Laby$ 
are translated into deviations of suitable transportation maps $\sigma_{xy}\colon A^Y\to A^X$ (mapping 
$e_y\mapsto e_x$ and $e_z\mapsto 0$ for $z\neq y$). So, for instance,  
\begin{multline*}
\xymatrixcolsep{1pc} \xymatrixrowsep{1pc}
\Phi(F)
\begin{bmatrix}  \xymatrix{
1 & \ar[l]_(0.7){a} \ar[dl]_(0.7){b}  1 \\
2 &  \ar[l]^(0.7){c} 2
} \end{bmatrix} 
 = F(a\sigma_{11}\de b\sigma_{21}\de c\sigma_{22}) \\ 
= F(a\sigma_{11}+b\sigma_{21}+ c\sigma_{22}) - F(a\sigma_{11}+ b\sigma_{21}) 
- F(a\sigma_{11}+ c\sigma_{22}) 
- F( b\sigma_{21}+ c\sigma_{22}) \\ 
+ F(a\sigma_{11}) +  F( b\sigma_{21}) + F(c\sigma_{22}) - F(0); \qquad a,b,c\in A.
\end{multline*}
Once again, polynomial functors of degree $n$ are those annihilated on sets of cardinality surpassing $n$ 
(\emph{ibid.}, Theorem 13).

The case of functors $\FAb\to\Ab$ is recovered by labelling all passages by $1$ (\emph{ibid.}, Theorem 14), 
commissioning the perspective of viewing e.g. the maze
$$
\xymatrixcolsep{1pc} \xymatrixrowsep{1pc}
\begin{bmatrix}  \xymatrix{
1 & \ar[l] \ar[dl]  1 \\
2 &  \ar[l] 2
} \end{bmatrix}
$$
(all labels equal to $1$) as a \emph{span} of surjections (read from right to left) 
\beq				\label{E: Span}
\maze{\{1,2\}}{\{a,b,c\}}{\{1,2\}}{f}{g}{}, \qquad 
f\colon 
\begin{cases}
a &\mapsto 1 \\ 
b,c &\mapsto 2,
\end{cases}
\quad 
g\colon 
\begin{cases}
a,b &\mapsto 1 \\ 
c &\mapsto 2,
\end{cases}
\eeq
which is substantially the road taken in \cite{BDFP}. Equivalence of the two approaches is proven in Theorem 11 of \cite{Laby}.

In order to extend the frame-work to functors $F\colon C\to D$, with $C$ a pointed algebraic theory and $D$ 
abelian, a modified labyrinth category $\Laby(C)$ is called for. 
The objects of $\Laby(C)$ are still finite sets, but a span like the one 
in \eref{E: Span} above will now also incorporate an arrow 
$\alpha\colon\Omega_1+\Omega_2\to\Omega_a+\Omega_b+\Omega_c$ of $C$. We write this symbolically as:
\beq				\label{E: 12-maze}
P=\maze{\{1,2\}}{\{a,b,c\}}{\{1,2\}}{f}{g}{\alpha}.
\eeq
Naturally, $\alpha$ will be subject to certain conditions, given in Definition \ref{D: Laby} below.
We may now enunciate our main results.

\begin{inttheorem}[\ref{S: Fun}]
The functor $\Phi\colon \Fun(C,D)\to \Add(\Laby(C),D)$, with the additive functor 
$\Phi(F)\colon \Laby(C)\to D$ given by 
the formul\ae 
$$
X \mapsto \ce_X F(\Omega,\dots,\Omega), \qquad
\maze{X}{Y}{Z}{f}{g}{\alpha} \mapsto F\left(\De_{y\in Y} \sigma_{f(y)y}\right)  F(\alpha),
$$
provides an equivalence of categories. 
\end{inttheorem}

So, for example, assuming $P$ as in \eref{E: 12-maze}, we would have 
\begin{multline*}
\Phi(F)(P) = F(\sigma_{1a}\de \sigma_{2b}\de \sigma_{2c})F(\alpha) \\ 
= \Big( F(\sigma_{1a}+ \sigma_{2b}+ \sigma_{2c})
-F(\sigma_{1a}+ \sigma_{2b})
-F(\sigma_{1a}+  \sigma_{2c}) \\
-F(\sigma_{2b}+ \sigma_{2c}) 
+F(\sigma_{1a})
+F(\sigma_{2b})
+F(\sigma_{2c})
-F(0) \Big) F(\alpha).
\end{multline*}

\begin{inttheorem}[\ref{S: Pol_n}]
The functor $F\colon C\to D$ is polynomial of degree $n$ if and only if $\Phi(F)$ annihilates sets of cardinality exceeding $n$.
\end{inttheorem}

The two types of labyrinth category, $\Laby(A)$, for a ring $A$, and $\Laby(C)$, for a category $C$, are found to 
agree in case $C$ is additive:

\begin{inttheorem}[\ref{S: Laby}]
When $C$ is an additive category, then $\Laby(C)$ is equivalent to $\Laby(C(\Omega,\Omega))$.
In particular, $\Laby(\FMod_A)$ is equivalent to $\Laby(A)$.
\end{inttheorem}

The special case of quadratic functors was investigated by Hartl \& Vespa in \cite{HV}, 
and we acquire their classification as a corollary in Theorem \ref{S: Quad} below.
In addition to vastly extending their result to functors of any degree, the proof is materially 
shortened and simplified. 

Since the category ${}_P\mathfrak{FAlg}$, for $P$ a set-operad, considered in \cite{HPV}, 
possesses finite sums (the free functor is left adjoint to the forgetful functor from the 
category of sets, which means it preserves sums), our main theorem applies to it. In principle, the results of \cite{HPV} 
could be derived, but we have chosen not to pursue this direction. 

Finally, let $B$ be a category with finite sums and a zero object, 
possessing a small, regular projective generator $\Omega$, along with all sums of it. Letting 
$C$ denote the subcategory of finite sums of copies of $\Omega$, 
one may consider the problem of extending 
functors defined on $C$ to the whole of $B$. 

\begin{inttheorem}[\ref{S: Fun Ext}]
A functor $F\colon C\to D$ has a unique (up to isomorphism) 
extension to a functor $\hat F\colon B\to D$ that is right-exact and preserves filtered inductive limits.
\end{inttheorem}

Uniqueness was known to Hartl \& Vespa (\cite{HV}, Proposition 6.25).  

We are greatly indebted to the Mittag--Leffler Institute for financial support. 
We would furthermore like to thank Djalal Mirmohades for assistance with the 
abstract nonsense.

\section{Pointed Algebraic Theories}

Throughout this paper, $C$ will denote a \textbf{pointed algebraic theory}. By this we mean a category with finite sums $+$  
and a zero object $0$, each object being isomorphic with a finite sum of copies of a generating 
object $\Omega$.

For any object $M\in C$, there are canonical zero arrows $0\to M$ and $M\to 0$. We shall denote all of these by simply $0$.
The identity arrow on $M$ will be written $1_M$.
Given objects $M_1,\dots, M_k$, there are canonical \emph{injections} $\iota_{M_i}\colon M_i\to M_1+\cdots+M_k$. 
We shall sometimes write simply $\iota_i$ if no confusion can arise.

Any finite family of arrows $\alpha_i\colon M_i\to N$, for $i=1,\dots,k$, possesses a \emph{sum} 
$$
\alpha_1+\cdots+\alpha_k \colon M_1+\cdots+M_k\to N.
$$
This way one defines the \emph{retractions} 
$$
\rho_i=\rho_{M_i} \colon M_1+\cdots+M_k\to M_i, \qquad 
\rho_i\iota_j = 
\begin{cases}
1_{M_i} & \text{if $i=j$;} \\
0 & \text{if $i\neq j$.}
\end{cases}
$$

Moreover, given a family of arrows $\alpha_i\colon M_i\to N_i$, for $i=1,\dots,k$, we define their ``direct sum'' 
$$
\alpha_1\oplus\cdots\oplus\alpha_k = \iota_{N_1}\alpha_1 + \cdots + \iota_{N_k}\alpha_k \colon M_1+\cdots +M_k \to N_1+\cdots+N_k.
$$ 

Finally, given two indexing sets $X$, $Y$, elements $x\in X$ and $y\in Y$, along with an object $M$, we define 
the \emph{transportation}
$$
\sigma_{xy}=\iota_x\rho_y\colon M^Y \to M \to M^X.
$$

\section{Deviations and Cross-Effects}				\label{A: Cross-Effects}

We develop, in this section, some general theory concerning 
functors $F\colon C\to D$ into an abelian category $D$. 
Write $[k]$ for the set $\{1,\dots,k\}$.

\bdf
The \textbf{deviation} of the arrows $\alpha_i\colon M_i\to N$ ($i=1,\dots, k$) in $C$ is 
$$
F(\alpha_1\de\cdots\de\alpha_k) = \sum_{I\subseteq [k]} (-1)^{k-\abs{I}} F\left(\sum_{i\in I} \alpha_i\right) \colon 
F(M_1+\cdots+M_k) \to F(N).
$$
\edf 

In contrast to the situation for e.g.~maps of modules (cf.~\cite{PolyMaps}), in which the deviation of any family 
of homomorphisms $M\to N$ exists, we only define the deviation in the very restricted case given above. 
A more involved approach in terms of a generalised Passi functor is possible, executed in Section 5 of \cite{HPV}.

Pivotal to our \emph{modus operandi} is the next formula. 
In case of maps of modules, a more general, yet simpler, formula is known (\cite{Laby}, Theorem 2).

We write $K\sqsubseteq A\times B$ to indicate that the set $K\subseteq A\times B$ 
and that the natural projections on $A$ and $B$ are onto.

\bth[The Deviation Formula]
Let $\alpha_i\colon N_i\to M$ ($i=1,\dots, m$) and $\beta_j\colon P_j\to \bigoplus_{i=1}^m N_i$ ($j=1,\dots, n$) be two families 
of arrows in $C$. Then 
\beq					\label{E: Dev Formula}
F(\alpha_1\de\cdots\de\alpha_m) \circ F( \beta_1\de\cdots\de\beta_n) 
= \sum_{K\sqsubseteq [m]\times [n]} \sum_{L\subseteq K} (-1)^{\abs{K}-\abs{L}}
 F\left( \sum_{j=1}^n \left( \sum_{(i,j)\in L} \alpha_i\right) \beta_j \right).
\eeq
\eth

\bpr
The right-hand member of \eref{E: Dev Formula} equals 
$$
\sum_{L\subseteq [m]\times [n]} (-1)^{\abs{L}}
 F\left( \sum_{j=1}^n \left( \sum_{(i,j)\in L} \alpha_i\right) \beta_j \right) 
 \left( \sum_{L\subseteq K\sqsubseteq [m]\times[n]} (-1)^{\abs{K}}\right).
$$
By the purely combinatorial Lemmata 1 and 2 of \cite{Laby},   
the inner sum reduces to $(-1)^{m+n+\abs{I}+\abs{J}+\abs{I}\abs{J}}$ when $L$ is of the form $I\times J$; 
else it is zero. The expression simplifies to 
$$
\sum_{I\times J\subseteq [m]\times [n]} (-1)^{m+n+\abs{I}+\abs{J}}
 F\left( \sum_{i\in I} \alpha_i \sum_{j\in J}\beta_j \right),
$$
which is the left-hand side of \eref{E: Dev Formula}.
\epr

The cross-effects of $F\colon C\to D$, originally defined by Eilenberg \& Mac Lane \cite{EM} for abelian $C$, 
were generalised to algebraic theories by Hartl \& Vespa \cite{HV}. As per their Proposition 1.3, they may be defined as 
follows. 

\bdf 
For objects $M_1,\dots,M_k\in C$, consider the canonical retractions $\hat\rho_j\colon \sum_{i} M_i \to \sum_{i\neq j} M_i$.
The \textbf{cross-effect} $\ce_k F(M_1,\dots,M_k)$ is the kernel of
$$ 
 (F(\hat\rho_1),\dots,F(\hat\rho_k)) \colon F(M_1+\cdots+ M_k) \to \prod_{j=1}^k F\left(\sum_{i\neq j} M_i\right).
$$
\edf

The definition of cross-effect under consideration can be brought to co-incide with that 
sometimes given for ordinary 
module functors, e.g. in \cite{PolyFun}.

\bth
The cross-effect $\ce_k F(M_1,\dots,M_k)$ equals the image of 
$$
F(\iota_1 \de \cdots \de \iota_k) \colon F(M_1+\cdots+ M_k) \to F(M_1+\cdots+ M_k).
$$
\eth 

\bpr
One easily verifies that
$$
F(\hat\rho_j)\circ F(\iota_1 \de \cdots \de \iota_k) = 
F(\hat\rho_j)\circ \sum_{I\subseteq [k]} (-1)^{k-\abs{I}} F\left(\sum_{i\in I} \iota_i\right) = 0, 
$$
since terms corresponding to $I$ with and without $j$ will cancel in the sum. This leads to a factorisation of 
$F(\iota_1 \de \cdots \de \iota_k)$ through the kernel of $(F(\hat\rho_1),\dots,F(\hat\rho_k))$: 
$$
\xymatrix{ 
F(M_1+\cdots+ M_k) \ar[dr]_{\alpha} \cruvear[ddr]_{\beta} 
\ar[rr]^{F(\iota_1 \de \cdots \de \iota_k)} && F(M_1+\cdots+ M_k) \\ 
& \Ker (F(\hat\rho_1),\dots,F(\hat\rho_k)) \ar[ur]_{\kappa} \ar@{.>}[d]^{\xi} \\ 
& L \cruvear[uur]_{\lambda}
}
$$
In the diagram, $\kappa$ is monic, $F(\iota_1 \de \cdots \de \iota_k)=\kappa\alpha$, and 
$$
F(\iota_1 \de \cdots \de \iota_k)\circ \kappa = 
\sum_{I\subseteq [k]} (-1)^{k-\abs{I}} F\left(\sum_{i\in I} \iota_i\right) \circ \kappa = 
F\left(\sum_{i\in[k]} \iota_i\right) \circ \kappa = \kappa
$$
because $\sum_{i\in I} \iota_i$, for $I\subset [k]$, factorises through $\hat\rho_j$ for any $j\notin I$.

In order to establish that $\Ker (F(\hat\rho_1),\dots,F(\hat\rho_k)) = \Im F(\iota_1 \de \cdots \de \iota_k)$, 
suppose $\lambda$ is another monic arrow such that $F(\iota_1 \de \cdots \de \iota_k)=\lambda\beta$. 
Then, since  $\kappa=F(\iota_1 \de \cdots \de \iota_k)\kappa=\lambda\beta\kappa$, the arrow $\xi$ fulfils 
$\kappa=\lambda\xi$ if and only if $\xi=\beta\kappa$. The existence of a unique factorisation settles the claim.  
\epr

The cross-effect $\ce_k F$ of rank $k$ is a multi-functor $C^k\to D$, 
and the functor $\ce_k\colon \Fun(C,D)\to\Fun(C^k,D)$ is exact (\cite{HV}, Proposition 1.5). 
As in the abelian case, there is a \emph{Cross-Effect Decomposition} (\emph{ibid.}, Proposition 1.4):
$$
F(M_1+\dots+ M_k) \cong \bigoplus_{I\subseteq [k]} \ce_I F((M_i)_{i\in I}).
$$
The idempotent corresponding to  $\ce_I F((M_i)_{i\in I})$ is $F(\De_{i\in I} \iota_i)$. 

\bdf
The functor $F$ is \textbf{polynomial} of degree $n$ if  $\ce_{n+1} F = 0$. 
Write $\Pol_n(C,D)$ for the full subcategory of $\Fun(C,D)$ of polynomial functors.
\edf

The category $\Pol_n(C,D)$ is closed under subfunctors, quotients and extensions (\cite{HV}, Proposition 1.7).

\section{The Labyrinth Category}				

Let $X$ and $Y$ be finite sets. The symbol 
$\con{X}{Y}{f}$,
with $f\colon Y\to X$ being a surjection, will be called a \textbf{contraction} from $Y$ to $X$. 
An \textbf{extension} from $Y$ to $X$ [\emph{sic!}] is a symbol of type 
$\ext{X}{Y}{g}{\alpha}$,
with $g\colon X\to Y$ being a surjection and $\alpha\colon \Omega^Y\to\Omega^X$ a
direct sum of arrows $\alpha\iota_y\colon \Omega\to\Omega^{g^{-1}(y)}$, for $y\in Y$.
The arrow $\alpha$ will be called the \textbf{structure map} of the extension. In what follows, 
it will sometimes be convenient to specify such a structure map as 
$\alpha\colon\Omega^{\tilde Y}\to \Omega^{\tilde X}$, for some supersets $\tilde Y\supseteq Y$, 
$\tilde X\supseteq X$. It is to be understood that the restriction to 
$\Omega^Y\to\Omega^X$ is all that matters. 

We shall consider the free additive semi-category having (formal) direct sums of finite sets for objects, 
and generated by extensions and contractions formed in accordance with the above rules. 
Certain axioms are now imposed.


\brom

\item $\con{X}{Y}{f}\con{Y}{Z}{g} = \con{X}{Z}{fg}$.

\item $\ext{X}{Y}{f}{\alpha}\ext{Y}{Z}{g}{\beta} = \ext{X}{Z}{gf}{\alpha\beta}$.

\item $\displaystyle \ext{X}{Y}{f}{\alpha}\con{Y}{Z}{g} = 
\sum_{K\sqsubseteq X\times_Y Z} \con{X}{K}{}   \lext{K}{Z}{}{\bigoplus_z \alpha\iota_{g(z)}}$,
where 
$$
\alpha\iota_{g(z)}\colon \Omega_z\to\Omega^Y\to\Omega^X, \qquad 
\bigoplus_z \alpha\iota_{g(z)} \colon \Omega^Z \to \Omega^{X\times Z},
$$ 
and
$$
X \times_Y Z = \set{(x,z)\in X\times Z | f(x) = g(z) }
$$ 
denotes the pull-back in the category of sets. 

\item $\con{X}{X}{1} = \maze{X}{X}{X}{1}{1}{1} = \ext{X}{X}{1}{1}$.

\erom

By Axioms \textsc{i}--\textsc{iv}, a spanning set for $\Laby$ consists of 
quantities of type \con{X}{Y}{f}\ext{Y}{Z}{g}{\alpha}, 
which we agree to write as simply $\maze{X}{Y}{Z}{f}{g}{\alpha}$.
The element exhibited in Axiom \textsc{iv} 
will function as the identical arrow at $X$, so we actually have a categorical structure. 

\brom

\setcounter{enumi}{4}

\item  We make the identification 
$$
\maze{X}{Y}{Z}{f}{g}{\alpha} = \maze{X}{\tilde Y}{Z}{\tilde f}{\tilde g}{\tilde \alpha},
$$
provided there be a bijection $Y\to \tilde Y$, making the diagrams commute: 
$$
\xymatrixcolsep{2pc}
\xymatrix{ 
X \ar[d]_{1} & \ar[l]_{f} Y \ar[r]^{g} \ar[d] & Z \ar[d]_{1} 
&& \Omega^Y \ar[d] & \ar[l]_{\alpha} \Omega^Z \ar[d]^{1} \\
X & \ar[l]_{\tilde f} \tilde Y \ar[r]^{\tilde g} & Z  
&& \Omega^{\tilde Y} & \ar[l]_{\tilde\alpha} \Omega^Z
}
$$
%
%

\item $\ext{X}{Y}{f}{\alpha}=0$ whenever there is a factorisation $\alpha\colon\Omega^Y\to\Omega^{\tilde X}\to\Omega^X$ 
for some proper subset $\tilde X\subset X$.

\item 
The equality
\begin{multline*}
\maze{P\cup\{p\}}{P\cup\{p\}}{Z}{}{g}{\beta} = 
\maze{P\cup\{p\}}{P\cup\{p_1,p_2\}}{Z}{}{\tilde g}{\alpha} \\
+ \maze{P\cup\{p\}}{P\cup\{p_1\}}{Z}{}{\tilde g}{\alpha}  
+ \maze{P\cup\{p\}}{P\cup\{p_2\}}{Z}{}{\tilde g}{\alpha}
\end{multline*}
is postulated to hold, provided that $p,p_1,p_2\notin P$, 
that the contracting surjections be the canonical ones, and that there be factorisations: 
$$
\xymatrix{ 
\Omega^Z \ar[r]^-{\alpha} \ar[dr]_-{\beta} & \Omega^P+\Omega_{p_1}+\Omega_{p_2} 
\ar[d]^{1_{\Omega^P}\oplus (\sigma_{pp_1}+\sigma_{pp_2})} 
&&&& P\cup\{p_1,p_2\}\ar[r] \ar[r] \ar[dr]_-{\tilde g} &  P\cup\{p\} \ar[d]^{g} \\
& \Omega^P + \Omega_p 
&&&&& Z
}
$$

\erom

\bdf				\label{D: Laby}
We call the category thus defined the \textbf{Labyrinth Category} over $C$, denoted by $\Laby(C)$, or simply $\Laby$.
\edf 

By construction, it is clear that $\Laby$ will be an additive category.
Axioms \textsc{vi} and \textsc{vii} fing their direct analogues in the two axioms for the 
Labyrinth Category $\Laby(A)$ founded upon a base ring $A$; cf.~\cite{Laby}, Definition 4.

The following extension of Axiom \textsc{vi} shall be needed shortly.

\blem			\label{L: GenAx}
Let $X$ and $Z$ be sets and $(Y_s)_{s\in S}$ a family of sets. Let $f\colon S\to X$ and $g\colon S\to Z$ be surjections. 
Then 
$$
\maze{X}{S}{Z}{f}{g}{\beta} = \sum_{(\emptyset\subset I_s\subseteq Y_s)_{s\in S}} 
\maze{X}{\set{(s,i) | s\in S,\ i\in I_s}}{Z}{\tilde f}{\tilde g}{\alpha} ,
$$
where
\begin{itemize} 
\item there is a factorisation 
$$ 
\beta\colon \xymatrix{\Omega^Z\ar[r]^-{\alpha} & \displaystyle \sum_{s\in S}\Omega^{Y_s} \ar[r] & \displaystyle \sum_{s\in S} \Omega}
$$
through a direct sum of folding maps $1+\cdots+1\colon \Omega^{Y_s}\to\Omega$;
\item the surjections $\tilde f$ and $\tilde g$ are projection on $S$ post-composed with $f$ and $g$, respectively. 
\end{itemize}
\elem

\bpr
Axiom \textsc{vi} may be inductively extended to 
$$
\maze{P\cup\{p\}}{P\cup\{p\}}{Z}{}{g}{\beta} = \sum_{\emptyset\subset I\subseteq [k]} 
\maze{P\cup\{p\}}{P\cup\set{p_i | i\in I}}{Z}{}{\tilde g}{\alpha} ,
$$
given factorisations: 
$$
\xymatrix{ 
\Omega^Z \ar[r]^-{\alpha} \ar[dr]_-{\beta} & \Omega^P+\Omega_{p_1}+\cdots+\Omega_{p_k} 
\ar[d]^{1_{\Omega^P}\oplus (\sigma_{pp_1}+\cdots+\sigma_{pp_k})} 
&&& P\cup\{p_1,\dots,p_k\}\ar[r] \ar[r] \ar[dr]_-{\tilde g} &  P\cup\{p\} \ar[d]^{g} \\
& \Omega^P + \Omega_p 
&&&& Z
}
$$
Post-composition with $\con{X}{P\cup\{p\}}{f}$ produces the identity
\beq					\label{E: Equation}
\maze{X}{P\cup\{p\}}{Z}{f}{g}{\beta} = \sum_{\emptyset\subset I\subseteq [k]} 
\maze{X}{P\cup\set{p_i | i\in I}}{Z}{\tilde f}{\tilde g}{\alpha},
\eeq
where 
$$
\tilde f \colon \xymatrix{P\cup \{p_1,\dots,p_k\} \ar[r] & P\cup\{p\} \ar[r]^-f & X }
$$
is the obvious map. Repeated application of \eref{E: Equation} yields the lemma. 
\epr

\section{General Functors}

For abelian $D$, we now propose to define a category equivalence 
$$
\Phi\colon \Fun(C,D)\to \Add(\Laby,D)
$$ 
between general functors on $C$ and 
additive functors on $\Laby(C)$. 

As a first step in the construction, we tentatively define, for a functor $F\colon C\to D$, an additive functor
$\Phi(F)\colon \Laby\to D$ by the equations
\begin{gather}
X \mapsto \ce_X F(\Omega,\dots,\Omega) \label{E: H(X)} \\ 
\con{X}{Y}{f} \mapsto F\left(\De_{y\in Y} \sigma_{f(y)y}\right) \\ 
\ext{Y}{Z}{g}{\alpha} \mapsto F\left(\De_{y\in Y} \sigma_{yy}\right) F(\alpha) = F\left(\De_{y\in Y} \sigma_{yy}\right) F(\alpha) 
F\left(\De_{z\in Z} \sigma_{zz}\right). \label{E: H(ext)}
\end{gather}
To perceive the equality of the two expressions in formula \eref{E: H(ext)}, note that 
\begin{multline*}
F\left(\De_{y\in Y} \sigma_{yy}\right) F(\alpha) F\left(\De_{z\in Z} \sigma_{zz}\right) \\
= \sum_{I\subseteq Y} (-1)^{\abs{Y}-\abs{I}} \sum_{J\subseteq Z} (-1)^{\abs{Z}-\abs{J}} 
F\left(\sum_{y\in I} \sigma_{yy}\right) F(\alpha) F\left(\sum_{z\in J} \sigma_{zz}\right) \\
= \sum_{I\subseteq Y} (-1)^{\abs{Y}-\abs{I}} 
F\left(\sum_{y\in I} \sigma_{yy}\right) F(\alpha) F\left(\sum_{z\in Z} \sigma_{zz}\right) 
=  F\left(\De_{y\in Y} \sigma_{yy}\right) F(\alpha),
\end{multline*}
since if $\hat z\in Z\setminus J$, then terms arising from $I$ and $I\cup\{\hat y\}$, 
where $\hat y\in Y$ is some element 
with $g(\hat y)=\hat z$, will cancel in the sum.

\blem 
Precisely one functor $\Phi(F)$ satisfies conditions \eref{E: H(X)}--\eref{E: H(ext)} above, and 
$$
\Phi(F)\maze{X}{Y}{Z}{f}{g}{\alpha} = F\left(\De_{y\in Y} \sigma_{f(y)y}\right)  F(\alpha).
$$
\elem

\bpr 
The formula is clear from the definition of $\Phi(F)$, provided it be well defined. We verify that $\Phi(F)$ respects the axioms 
of $\Laby$. 

\emph{Axiom I.} 
By the Deviation Formula, 
\begin{multline*}
\Phi(F)\con{X}{Y}{f} \circ \Phi(F) \con{Y}{Z}{g} = F\left(\De_{y} \sigma_{f(y)y}\right) \circ F\left(\De_{z} \sigma_{g(z)z}\right) \\
= \sum_{K\sqsubseteq Y\times Z} \sum_{L\subseteq K} (-1)^{\abs{K}-\abs{L}} 
F\left( \sum_{z\in Z} \left(\sum_{(y,z)\in L} \sigma_{f(y)y}\right) \sigma_{g(z)z}\right).
\end{multline*}
If $K$ contains some $(\hat y,\hat z)$ with $g(\hat z)\neq \hat y$, then terms corresponding to $L$ with and without 
$(\hat y,\hat z)$ cancel. Therefore only terms with $K=\set{(g(z),z)|z\in Z}$ survive, and the expression simplifies to 
\begin{multline*}
\sum_{L\subseteq \set{(g(z),z)|\, z\in Z}} (-1)^{\abs{Z}-\abs{L}} 
F\left( \sum_{z\in Z} \left(\sum_{(g(z),z)\in L} \sigma_{fg(z)g(z)}\right) \sigma_{g(z)z}\right) \\
= \sum_{L\subseteq Z} (-1)^{\abs{Z}-\abs{L}} 
F\left( \sum_{z\in L} \sigma_{fg(z)z}\right) = \Phi(F)\con{X}{Z}{fg}.
\end{multline*}

\emph{Axiom II.} 
Using the equality in formula \eref{E: H(ext)} above, we find
\begin{multline*}
\Phi(F)\ext{X}{Y}{f}{\alpha} \circ \Phi(F) \ext{Y}{Z}{g}{\beta} 
= F\left(\De_{x} \sigma_{xx}\right) F(\alpha) \circ F\left(\De_{y} \sigma_{yy}\right) F(\beta) \\
= F\left(\De_{x} \sigma_{xx}\right) F(\alpha) F(\beta) = \Phi(F)\ext{X}{Z}{gf}{\alpha\beta}.
\end{multline*}

\emph{Axiom III.} 
One easily verifies that 
\begin{multline*}
\Phi(F)\ext{X}{Y}{f}{\alpha} \circ \Phi(F)\con{Y}{Z}{g} = 
F\left(\De_{x} \sigma_{xx}\right) F(\alpha) \circ F\left(\De_{z} \sigma_{g(z)z}\right) \\
= F\left(\De_{x} \sigma_{xx}\right) F\left(\De_z \sum_{(x,z)\in X\times_Y Z} \sigma_{x,(x,z)}\right) 
F\left(\bigoplus_z \alpha\iota_{g(z)}\right) ,
\end{multline*}
which, by the Deviation Formula, equals
$$
 \sum_{K\sqsubseteq X\times Z} \sum_{L\subseteq K} (-1)^{\abs{K}-\abs{L}} 
F\left( \sum_{z} \left( \sum_{(\tilde x,z)\in L} 
\sigma_{\tilde x,\tilde x} \sum_{(x,z)\in X\times_Y Z} \sigma_{x,(x,z)} \right) \right) 
F\left(\bigoplus_z \alpha\iota_{g(z)}\right). 
$$
If $K$ contains a pair $(\hat x,\hat z)\notin X\times_Y Z$, terms corresponding to $L$ with and 
without $(\hat x,\hat z)$ 
will cancel. Thus $K\sqsubseteq X\times_Y Z$, and the expression reduces to 
\begin{multline*}
\sum_{K\sqsubseteq X\times_Y Z} \sum_{L\subseteq K} (-1)^{\abs{K}-\abs{L}} 
F\left( \sum_{(x,z)\in L}  \sigma_{x,(x,z)} \right) F\left(\bigoplus_z \alpha\iota_{g(z)}\right) = \\
\sum_{K\sqsubseteq X\times_Y Z} 
F\left( \De_{(x,z)\in K}  \sigma_{x,(x,z)} \right) F\left(\bigoplus_z \alpha\iota_{g(z)}\right) 
= \sum_{K\sqsubseteq X\times_Y Z} \Phi(F)\lmaze{X}{K}{Z}{}{}{\bigoplus_z \alpha\iota_{g(z)}}.
\end{multline*}

\emph{Axiom IV, V.} Clear.

\emph{Axiom VI.}  
Supposing there is a factorisation $\alpha\colon\Omega^Y\to\Omega^{\tilde X}\to\Omega^X$ 
for some proper subset $\tilde X\subset X$, let $\hat x\in X\setminus\tilde X$. Then
$$
\Phi(F)\ext{X}{Y}{f}{\alpha} = F\left(\De_{x} \sigma_{xx}\right) F(\alpha) = 
\sum_{I\subseteq X} (-1)^{\abs{X}-\abs{I}} F\left(\sum_{x\in I} \sigma_{xx}\right) F(\alpha) = 0,
$$
since terms arising from $I$ with and without $\hat x$ will cancel in the sum.

\emph{Axiom VII.} 
Assuming ourselves placed in the position described above for Axiom \textsc{vii}, we have 
\begin{multline*}
\Phi(F)\maze{P\cup\{p\}}{P\cup\{p_1,p_2\}}{Z}{}{\tilde g}{\alpha} 
+ \Phi(F)\maze{P\cup\{p\}}{P\cup\{p_1\}}{Z}{}{\tilde g}{\alpha}  \\
+ \Phi(F)\maze{P\cup\{p\}}{P\cup\{p_2\}}{Z}{}{\tilde g}{\alpha} =  \\
 F\left(\De_{s\in P}\sigma_{ss}\de\sigma_{pp_1}\de\sigma_{pp_2}\right) F(\alpha)  + 
F\left( \De_{s\in P}\sigma_{ss}\de\sigma_{pp_1} \right) F(\alpha)  + 
F\left( \De_{s\in P}\sigma_{ss}\de\sigma_{pp_2} \right) F(\alpha) =\\
 F\left( \De_{s\in P}\sigma_{ss}\de\sigma_{pp} \right) F(1_{\Omega^P}\oplus(\sigma_{pp_1}+\sigma_{pp_2})) F(\alpha) =
\Phi(F)\maze{P\cup\{p\}}{P\cup\{p\}}{Z}{}{g}{\beta}.
\end{multline*}
This concludes the proof. 
\epr

We next define the action of $\Phi$ on natural transformations. Given $\eta\colon F\to G$, let 
$\Phi(\eta)\colon\Phi(F)\to\Phi(G)$ be the restriction of the appropriate cross-effect: 
$$
\Phi(\eta)_X = \ce_X \eta \colon \ce_X F \to \ce_X G.
$$

\blem
 $\Phi\colon \Fun(C,D) \to \Add(\Laby,D)$ is a fully faithful functor. 
\elem

\bpr
A natural transformation $\eta$ 
can be unambiguously re-assembled from its components $\ce_X\eta$ because of the Cross-Effect Decomposition.
\epr

There remains to verify that $\Phi$ is essentially surjective. For an additive functor $H\colon \Laby\to D$, 
we define its pre-image $F\colon C\to D$ under $\Phi$ by 
$$
\Omega^X \mapsto \bigoplus_{A\subseteq X} H(A), \qquad 
\left[ \alpha\colon \Omega^Y \to \Omega^X \right] \mapsto 
\sum_{U\subseteq X\times Y} H\lmaze{U_X}{U}{U_Y}{}{}{\bigoplus_{y\in Y}\alpha\iota_y},
$$
where $U_X$ and $U_Y$ denote the projections of $U\subseteq X\times Y$ on the components $X$ and $Y$, respectively. 

\blem
$F$ is a functor.  
\elem 

\bpr 
Consider two arrows $\alpha\colon\Omega^Y\to\Omega^X$ and $\beta\colon\Omega^Z\to\Omega^Y$ of $C$. We have 
\begin{multline*}
F(\alpha)F(\beta) = \sum_{U\subseteq X\times Y} H\lmaze{U_X}{U}{U_Y}{}{}{\bigoplus_{y\in Y}\alpha\iota_y} 
\circ \sum_{V\subseteq Y\times Z} H\lmaze{V_Y}{V}{V_Z}{}{}{\bigoplus_{z\in Z}\beta\iota_z} \\ 
= \sum_{\substack{A\subseteq X \\ B\subseteq Y \\ E\subseteq Z}} 
\sum_{\substack{U\sqsubseteq A\times B \\ V\sqsubseteq B\times E}} 
H\left( \lmaze{A}{U}{B}{}{}{\bigoplus_{y\in Y}\alpha\iota_y} \lmaze{B}{V}{E}{}{}{\bigoplus_{z\in Z}\beta\iota_z} \right) \\ 
= \sum_{\substack{A\subseteq X \\ B\subseteq Y \\ E\subseteq Z}} 
\sum_{\substack{U\sqsubseteq A\times B \\ V\sqsubseteq B\times E}} 
\sum_{W\sqsubseteq U\times_B V} 
H\left( \con{A}{U}{} \llmaze{U}{W}{V}{}{}{\bigoplus\limits_{(y,z)\in Y\times Z}\alpha\iota_y} 
\lext{V}{E}{}{\bigoplus_{z\in Z}\beta\iota_z} \right) ,
\end{multline*}
using just the definition of $F$ and the rules for composition in $\Laby$. Now rearrange the sum and apply 
Lemma \ref{L: GenAx} (p.~\pageref{L: GenAx}):
\begin{multline*}
F(\alpha)F(\beta) = \sum_{W\subseteq X\times Y\times Z} 
H \llmaze{W_X}{W}{W_Z}{}{}{\bigoplus\limits_{z\in Z} \left(\bigoplus\limits_{y\in Y}\alpha\iota_y\right)\beta\iota_z}  \\ 
= \sum_{S\subseteq X\times Z} \sum_{\big(\emptyset\subset I_{xz}\subseteq Y\big)_{(x,z)\in S} }
H \llmaze{S_X}{\bigcup\limits_{(x,z)\in S} \Big(\{x\}\times I_{xz} \times \{z\}\Big)}{S_Z}{}{}{\bigoplus\limits_{z\in Z} 
\left(\bigoplus\limits_{y\in Y}\alpha\iota_y\right) \beta\iota_z}  \\
= \sum_{S\subseteq X\times Z} 
H \llmaze{S_X}{S}{S_Z}{}{}{\bigoplus_{z\in Z} \alpha\beta\iota_z} = F(\alpha\beta). \qedhere
\end{multline*}
\epr

\blem
$\Phi(F)=H$.
\elem

\bpr 
We calculate, using the definitions of $\Phi(F)$ and $F$:
\begin{multline*}
\Phi(F)\con{X}{Y}{f} = F\left(\De_{y\in Y} \sigma_{f(y)y}\right) = 
\sum_{I\subseteq Y} (-1)^{\abs{Y}-\abs{I}} F\left(\sum_{y\in I} \sigma_{f(y)y}\right) \\ 
= \sum_{I\subseteq Y} (-1)^{\abs{Y}-\abs{I}} \sum_{U\subseteq X\times I} H\llmaze{U_X}{U}{U_I}{}{}{\bigoplus_{y\in Y} 
\sigma_{f(y)y}} \\ 
= \sum_{U\subseteq X\times Y} H\llmaze{U_X}{U}{U_Y}{}{}{\bigoplus_{y\in Y} 
\sigma_{f(y)y}} \sum_{U_Y\subseteq I\subseteq Y} (-1)^{\abs{Y}-\abs{I}} .
\end{multline*}
The inner sum is $1$ if $U_Y=Y$ and $0$ otherwise. The structure map factorises
$$
\bigoplus_{y\in Y} \sigma_{f(y)y}\colon \Omega^Y \to \Omega^{\set{(f(y),y) | \, y\in Y}} \to \Omega^{X\times Y},
$$  
so, by Axiom \textsc{vi}, only the term $U=\set{(f(y),y) | y\in Y}$ survives, making 
$$
\Phi(F)\con{X}{Y}{f} = H\llmaze{X}{\set{(f(y),y) | y\in Y}}{Y}{}{}{\bigoplus_{y\in Y} \sigma_{f(y)y}} =  H\con{X}{Y}{f}.
$$

The formula 
$$
F\left(\De_{y\in Y} \sigma_{f(y)y}\right) = H\maze{X}{Y}{Y}{f}{}{},
$$ 
extracted from the above computation, contains, as a special case,  
$$
F\left(\De_{x\in X} \sigma_{xx}\right) = H\maze{X}{X}{X}{}{}{} ,
$$ 
which we now apply to consider
$$
\Phi(F)\ext{X}{Y}{f}{\alpha} = F\left(\De_{x\in X} \sigma_{xx}\right) F(\alpha) 
= \sum_{\substack{V\subseteq X\times Y \\ V_X=X}} H\lmaze{V_X}{V}{V_Y}{}{}{\bigoplus_{y\in Y} \alpha\iota_y}.
$$
Because $\alpha$ is a direct sum of arrows $\alpha\iota_y\colon \Omega\to\Omega^{f^{-1}(y)}$,
the structure map factorises
$$
\bigoplus_{y\in Y} \alpha\iota_y\colon \Omega^Y \to \Omega^{\set{(x,f(x))| \, x\in X}} \to \Omega^{X\times Y},
$$ 
so, by Axiom \textsc{vi}, only the term $V=\set{(x,f(x))|x\in X}$ survives, making
$$
\Phi(F)\ext{X}{Y}{f}{\alpha} = H\lmaze{X}{\set{(x,f(x))|x\in X}}{Y}{}{}{\bigoplus_{y\in Y} \alpha\iota_y} = H\ext{X}{Y}{f}{\alpha}.
$$
The proof is complete. 
\epr

Summarising our finds, we have thus proved our main theorem.

\bth					\label{S: Fun}
The functor $\Phi\colon \Fun(C,D)\to \Add(\Laby,D)$, with the additive functor $\Phi(F)\colon \Laby\to D$ given by 
the formul\ae 
$$
X \mapsto \ce_X F(\Omega,\dots,\Omega), \qquad
\maze{X}{Y}{Z}{f}{g}{\alpha} \mapsto F\left(\De_{y\in Y} \sigma_{f(y)y}\right)  F(\alpha),
$$
provides an equivalence of categories. 
\eth

\section{Polynomial Functors}

Polynomial functors find a natural interpretation in the context of Theorem~\ref{S: Fun}.

\bth				\label{S: Pol_n}
The functor $F\colon C\to D$ is polynomial of degree $n$ if and only if $\Phi(F)$ annihilates sets of cardinality exceeding $n$.
\eth

\bpr
Suppose $\Phi(F)(X)=0$ whenever $\abs{X}>n$. This means that $F(\iota_1\de\cdots\de\iota_k) = 0$, where
$\iota_i$ ($i=1,\dots, k$) are injections associated to a sum $\Omega^k$.
Consider any sum $M=M_1+\dots+ M_k$ of $k>n$ non-zero objects. Owing to the structure 
of the category $C$, each $M_i=\Omega^{m_i}$ for some positive integer $m_i$. 
Writing $\iota_{ij}\colon \Omega_j\to M_i\to M$ ($j=1,\dots, m_i$) for the natural injections,
we may compute 
\begin{multline*}
\ce_k F(M_1,\dots,M_k) =  \\ \Im F\left(\De_{i=1}^k \sum_{j=1}^{m_i} \iota_{ij}\right) 
= \Im \sum_{\emptyset\subset J_1\subseteq [m_1]} \!\cdots\! \sum_{\emptyset\subset J_k\subseteq [m_k]} 
F\left( \De_i \De_{j\in J_i} \iota_{ij} \right) = 0,
\end{multline*}
each inner term vanishing by the assumption placed on $\Phi(F)$. So $F$ is polynomial of degree $n$. 
The converse is clear from the definition of $\Phi(F)$. 
\epr

\section{The Case of Quadratic Functors}

The case of quadratic functors merits some attention. We denote by 
$\Laby_{(2)}=\Laby_{(2)}(C)$ the quotient of $\Laby(C)$ obtained by annihilating 
all sets of cardinality greater than $2$. Theorem \ref{S: Pol_n} provides a category equivalence 
$$
\Pol_2(C,D)\cong\Add(\Laby_{(2)},D).
$$ 
The category $\Laby_{(2)}$ has three isomorphism classes of objects: $\{\}$, $\{1\}$ and $\{1,2\}$. 
The empty set corresponds to the value of the functor $F\colon C\to D$ at the zero object. Assuming $F(0)=0$, this object may be 
disregarded. Between the sets $[1]=\{1\}$ and $[2]=\{1,2\}$, a spanning system may be compiled: 
\begin{align*}
I(\alpha,\beta) &= \lmaze{[2]}{[2]}{[2]}{}{}{\alpha\oplus\beta}, & \alpha,\beta\colon&\Omega\to\Omega \\ 
T(\gamma,\delta) &= \lmaze{[2]}{[2]}{[2]}{(12)}{}{\gamma\oplus\delta}, & \gamma,\delta\colon&\Omega\to\Omega \\
P(\zeta,\eta) &= \lmaze{[1]}{[2]}{[2]}{}{}{\zeta\oplus\eta}, & \zeta,\eta\colon&\Omega\to\Omega \\
H(\xi) &= \lmaze{[2]}{[2]}{[1]}{}{}{\xi}, & \xi\colon&\Omega\to\Omega^2  \\  
I(\epsilon) &= \lmaze{[1]}{[1]}{[1]}{}{}{\epsilon}, & \epsilon\colon&\Omega\to\Omega  \\  
E(\omega) &= \lmaze{[1]}{[2]}{[1]}{}{}{\omega}, & \omega\colon&\Omega\to\Omega^2  . 
\end{align*}
Valid laws of multiplication are as beneath: 
{\small 
\begin{align*}
I(\alpha,\beta)I(\tilde\alpha,\tilde\beta) &= I(\alpha\tilde\alpha,\beta\tilde\beta) 
& I(\epsilon)I(\tilde\epsilon) &= I(\epsilon\tilde\epsilon) \\
I(\alpha,\beta)T(\gamma,\delta) &= T(\beta\gamma,\alpha\delta) 
& I(\epsilon)E(\omega) &= E((\epsilon\oplus\epsilon)\omega) \\
T(\gamma,\delta)I(\alpha,\beta) &= T(\gamma\alpha,\delta\beta) 
& E(\omega)I(\epsilon) &= E(\omega\epsilon) \\
 T(\gamma,\delta)T(\tilde\gamma,\tilde\delta) &= I(\delta\tilde\gamma,\gamma\tilde\delta) 
& E(\omega)E(\tilde\omega) 
&=   E((\rho_1\omega\oplus\rho_2\omega)\tilde\omega) + E((\rho_2\omega\oplus\rho_1\omega)\tilde\omega) \\
I(\alpha,\beta)H(\xi) &= H((\alpha\oplus\beta)\xi) 
& H(\xi)E(\omega) &= H((\rho_1\xi\oplus\rho_2\xi)\omega) + H((\rho_1\xi\oplus\rho_2\xi)\sigma\omega) \\
T(\gamma,\delta)H(\xi) &= H((\delta\oplus\gamma)\xi) 
& H(\xi)P(\zeta,\eta) &= I(\rho_1\xi\zeta,\rho_2\xi\eta) + T(\rho_2\xi\zeta,\rho_1\xi\eta) \\
P(\zeta,\eta)I(\alpha,\beta) &= P(\zeta\alpha,\eta\beta) 
& E(\omega)P(\zeta,\eta) &= P(\rho_1\omega\zeta,\rho_2\omega\eta) + P(\rho_2\omega\zeta,\rho_1\xi\eta) \\
P(\zeta,\eta)T(\gamma,\delta) &= P(\eta\gamma,\zeta\delta) 
& P(\zeta,\eta)H(\xi) &= E((\zeta\oplus\eta)\xi) \\
H(\xi)I(\epsilon) &= H(\xi\epsilon) 
& I(\epsilon)P(\zeta,\eta) &= P(\epsilon\zeta,\epsilon\eta) ,
\end{align*}
}%
where we have abbreviated $\sigma=\sigma_{12}+\sigma_{21}$. 

Putting $T=T(1,1)$ and $P=P(1,1)$, we may express 
$$
T(\gamma,\delta)=TI(\gamma,\delta), \qquad P(\zeta,\eta)=PI(\zeta,\eta) \qquad\text{and}\qquad E(\omega) = PH(\omega). 
$$
As an additive category, $\Laby_{(2)}$ is generated by the mazes  $I(\epsilon)$, 
$I(\alpha,\beta)$, $H(\xi)$, $T$ and $P$, subject to the ten axioms: 
{\small 
\begin{align*}
I(\epsilon)I(\tilde\epsilon) &= I(\epsilon\tilde\epsilon) 								\Label{E: M module} &
I(\alpha,\beta)I(\tilde\alpha,\tilde\beta) &= I(\alpha\tilde\alpha,\beta\tilde\beta)		\Label{E: N module} \\
I(\alpha,\beta)T &= TI(\beta,\alpha) 													\Label{E: N symmetric} &
T^2&=I(1,1) 																			\Label{E: T involution} \\
I(\epsilon)P &= PI(\epsilon,\epsilon) 												\Label{E: P homo} &
PT&=P 																				\Label{E: P}  \\
H(\xi)I(\epsilon) &= H(\xi\epsilon) 													\Label{E: H well-defined} &
TH(\xi) &= H(\sigma\xi) 																\Label{E: H symmetric} \\
 I(\alpha,\beta)H(\xi) &= H((\alpha\oplus\beta)\xi) 									\Label{E: H homo} &
H(\xi)P &= I(\rho_1\xi,\rho_2\xi) + I(\rho_2\xi,\rho_1\xi) 							\Label{E: QM2} 
\end{align*}
}%

We thus retrieve Hartl \& Vespa's classification \cite{HV}, Theorem 7.1, of quadratic functors 
$F\colon C\to\Ab$.  

\bth				\label{S: Quad}
Let $T_{11}\ce_2 U$ denote the bilinearisation of the second cross-effect of the standard projective functor $U=\Z[C(\Omega,-)]$ 
and let $\Lambda$ be the reduced monoid ring of $C(\Omega,\Omega)$.
Quadratic functors $F\colon C\to \Ab$ satisfying $F(0)=0$ are equivalent to diagrams 
$$
\xymatrixcolsep{2pc}
\xymatrix{
T_{11}\ce_2 U(\Omega,\Omega) \otimes_{\Lambda} M_e \ar[r]^-{\hat H} & M_{ee} \ar@(ul,ur)^T \ar[r]^{P} & M_e  ,
}
$$
fulfilling the conditions: 
\begin{itemize}
\item $M_e$ is a left $\Lambda$-module. 
\item $M_{ee}$ is a symmetric left $\Lambda\otimes\Lambda$-module with involution $T$. 
\item $P$ is an homomorphism of $\Lambda$-modules (with respect to the diagonal action of $\Lambda$ on $M_{ee}$), 
satisfying $PT=P$. 
\item $\hat H$ is an homomorphism of symmetric $\Lambda\otimes\Lambda$-modules 
\end{itemize}
\eth

\bpr
Equations \eref{E: M module}--\eref{E: T involution} state that $\Phi(F)([1])=\ce_1 F(\Omega)$ is a module 
over $\Lambda$ and $\Phi(F)([2])=\ce_2 F(\Omega,\Omega)$ a symmetric module over $\Lambda\otimes\Lambda$ with involution~$T$. 

By equation \eref{E: P homo}, $P$ is an homomorphism of $\Lambda$-modules. 
Equations \eref{E: H well-defined}--\eref{E: H homo} express the fact that 
$$
\hat H\colon T_{11}\ce_2 U(\Omega,\Omega) \otimes_{\Lambda} M \to N, \qquad [\xi] \otimes x \mapsto H(\xi)x,
$$ 
is a well-defined homomorphism of symmetric $\Lambda\otimes\Lambda$-modules.
Equation \eref{E: QM2} is Hartl \& Vespa's condition (QM2). Their condition (QM1) is a recasting of the 
single pertinent instance of our Axiom \textsc{vii} for $\Laby$.
\epr

\section{The Case of an Additive Category}

The description of $\Laby(C)$ simplifies considerably in the special 
case of an additive category $C$, in which the sum of two objects is actually a biproduct. 
Such is the situation in the category $\FMod_A$ of finitely generated, free right modules over a ring $A$ 
(unital, but not necessarily commutative). 
Recall that the labyrinth category $\Laby(A)$ was originally defined for a ring $A$ in the paper~\cite{Laby}. 
This case turns out to be prototypical. 

\bth			\label{S: Laby}
When $C$ is an additive category, then $\Laby(C)$ is equivalent to $\Laby(C(\Omega,\Omega))$.
In particular, $\Laby(\FMod_A)$ is equivalent to $\Laby(A)$.
\eth

\bpr
One easily verifies the existence of a functor 
$$
\Xi\colon\Laby(C)\to\Laby(C(\Omega,\Omega)), \quad
\maze{X}{Y}{Z}{f}{g}{\alpha} \mapsto \bigcup_{y\in Y} 
\xymatrixcolsep{3pc} 
\{\xymatrix{f(y) & \ar[l]_-{\rho_{f(y)}\alpha\iota_{g(y)}} g(y) } \}.
$$ 
Let $\bigcup_{i=1}^k \{\xymatrix{x_i & \ar[l]_{\epsilon_i} z_i } \}$ be a maze with labels 
$\epsilon_i\colon \Omega_{z_i}\to \Omega_{x_i}$.
Denoting $X=\#\set{x_i|i\in [k]}$ and  $Z=\#\set{z_i|i\in [k]}$, where $\#A$ denotes the support of the multi-set $A$, 
the inverse of $\Xi$ transforms
$$
\bigcup_{i=1}^k \{\xymatrix{x_i & \ar[l]_{\epsilon_i} z_i } \} \mapsto
\maze{X}{[k]}{Z}{f}{g}{\alpha},
$$
where $f(i)=x_i$, $g(i)=z_i$, and the structure map $\alpha\colon \Omega^Z\to\Omega^k$ is given by 
$$
\xymatrixcolsep{2pc}
\alpha\iota_z \colon \xymatrix{
\Omega_z \ar[r]^-{\begin{pmatrix}1 \\ \vdots \\ 1 \end{pmatrix}} 
& \Omega^{\set{i \, | \, z_i=z}} \ar[r]^-{\bigoplus_i \epsilon_i} 
& \Omega^{\set{i \, | \, z_i=z}}, \qquad z\in Z. 
}
$$
Functoriality is readily checked and follows from first principles.
\epr

\section{Extension of Functors}

Let $B$ be a category with finite sums and a zero object. 
Writing $C$ for the subcategory of objects isomorphic with 
finite sums of a special object $\Omega\in B$, we shall now solve the problem of extending a functor $C\to D$ to 
a functor $B\to D$. As before, $D$ denotes an abelian category.  

Recalling that a \emph{regular epimorphism} is an epimorphism occurring as a co-equaliser, we make
 the following assumptions:

(a) $B$ has all (small) sums of $\Omega$. A sum $\Omega^\kappa$, for an ordinal $\kappa$, is called a 
\emph{free object}.

(b) $\Omega$ is a \emph{generator}, meaning that any object $M$ admits a regular epimorphism 
$P\to M$ from a free object. 

(c) $\Omega$ is \emph{regular projective}, in the sense that arrows out of $\Omega$ lift through  
regular epimorphisms.

(d) $\Omega$ is \emph{small}, in the sense that $B(\Omega,-)$ preserves filtered inductive limits. 

Some further terminology will be required, following \cite{HV}. A co-equaliser 
\beq					\label{E: Co-equaliser}
\xymatrix{ 
P \upar[r]^{\zeta} \downar[r]_{\eta} & Q \ar[r]^{\theta} & M,
}
\eeq
is \textbf{free} if both objects $P$ and $Q$ are. If any two arrows $\alpha,\beta\colon\Omega\to Q$ 
factorise as $\alpha=\zeta\gamma$ and $\beta=\eta\gamma$ for some $\gamma\colon\Omega\to P$, the 
co-equaliser is said to be \textbf{saturated}. 

Hartl \& Vespa make a convincing case (\cite{HV}, Section 6.5) for why a \textbf{right-exact} functor ought to be defined 
as a functor transforming saturated, free co-equalisers into co-equalisers.

%
%

\bth					\label{S: Fun Ext}
Under the assumptions (a)--(d) above, a functor $F\colon C\to D$ has a unique (up to isomorphism) 
extension to a functor $\hat F\colon B\to D$ that is right-exact and preserves filtered inductive limits.
\eth

\bpr
From the theory of categories, $F$ will extend uniquely to a functor 
$
\tilde F\colon \Fun(C^\op,\Sets)\to D$, 
preserving all inductive limits (\cite{Sheaves}, Corollary~I.5.4). Here $C$ embeds as usual in $\Fun(C^\op,\Sets)$  through 
the Yoneda embedding $\Upsilon\colon P\mapsto B(-,P)$. 
It has a canonical extension $J\colon B\to \Fun(C^\op,\Sets)$, given by the same formula $M\mapsto B(-,M)|_{C}$. 
$$
\xymatrixcolsep{3pc}
\xymatrix{ 
C \ar[d]_{F} \ar[r]^-{\Upsilon} & \Fun(C^\op,\Sets) \ar[dl]_{\tilde F} \\ 
D & B \ar[u]_{J} \ar[l]_{\hat F}
}
$$

We prove that $J$ is right-exact and preserves filtered inductive limits; the same must then hold true for the
composite $\hat F=\tilde F \circ J$. Limits are calculated point-wise, and in as much as $B(\Omega,-)$ 
(and therefore $B(\Omega^m,-)$ for any finite $m$) preserves filtered 
inductive limits (by (d) above), so will $J$. 

To prove $J$ preserves free, saturated co-equalisers, 
consider such a diagram together with an object $R\in C$: 
$$	
\xymatrixcolsep{2pc}		
\xymatrix{ 
& R \ar[dl]_{\gamma} \upar[d]^{\alpha} \downar[d]_{\beta} \ar[dr]^{\mu} 
&&&& K \\
P \upar[r]^{\zeta} \downar[r]_{\eta} & Q \ar[r]^{\theta} & M 
& B(R,P) \upar[r]^{\zeta_\ast} \downar[r]_{\eta_\ast} & B(R,Q) \ar[r]^{\theta_\ast} \ar[ur]^{\kappa} & B(R,M) \dotar[u]_{\lambda}
}
$$
Suppose $\kappa\zeta_\ast=\kappa\eta_\ast$,
and let $\mu\colon R\to M$. By (c) above, 
there will be a factorisation $\mu=\theta\alpha$, and, if the diagram is to commute, 
$\lambda(\mu)=\lambda(\theta\alpha)=\kappa(\alpha)$. 
If also $\mu=\theta\beta$, then, by saturation, $\alpha=\zeta\gamma$ and $\beta=\eta\gamma$ for some arrow $\gamma$. 
Now $\kappa(\alpha)=\kappa(\zeta\gamma)=\kappa(\eta\gamma)=\kappa(\beta)$, so there is a well-defined, unique 
$\lambda\colon \mu\mapsto\kappa(\alpha)$ making the diagram commute, establishing that $B(R,M)$ is the co-equaliser. 

The functor $\hat F$ extends $F$: for $P\in C$, one has 
$$
\hat F(P)=\tilde F J(P)=\tilde F(B(-,P))=\tilde F(\Upsilon(P)) = F(P).
$$

Uniqueness of $\hat F$ is clear from the following considerations. If $P=\Omega^{\kappa}$ is free, 
it is the filtered inductive limit of finite sums of $\Omega$, 
so we may write $P=\indlim_{i\in I} P_i$ with $P_i\in C$. Then, necessarily, 
$$
\hat F(P) = \hat F(\indlim P_i) = \indlim \hat F(P_i) = \indlim F(P_i),
$$ 
since $\hat F$ is required to commute with filtered inductive limits. Moreover, using property (b) above,
Hartl \& Vespa (\cite{HV}, Lemma 6.22) prove the existence of a free, saturated co-equaliser \eref{E: Co-equaliser} for 
an arbitrary object $M\in B$. Then 
$$
\hat F(M) = \Coker (\hat F(\zeta)-\hat F(\eta)),
$$
since $\hat F$ commutes with free, saturated co-equalisers, and the action of $\hat F$ on $\zeta$ and $\eta$ 
was just asserted to be determined by $F$. This concludes the proof.
\epr

\end{document}